\documentclass{amsart}

\usepackage{pgfplotstable}
\usepackage{booktabs}
\usepackage{array}
\usepackage{colortbl}
\usepackage{mathtools}
\usepackage{graphicx}
\usepackage{figsize}
\usepackage{enumerate}
\usepackage{amsmath,amssymb,latexsym, geometry}

\usepackage{color}
\usepackage{subfigure}
\usepackage{multirow}
\usepackage{mathptmx}
\usepackage{url,bm,overpic}
\usepackage{enumerate}
\usepackage{bbm}
\usepackage[ruled,vlined,linesnumbered]{algorithm2e}
\usepackage{comment}
\newtheorem{theorem}{Theorem}[section]

\newtheorem{proposition}{Proposition}[section]

\newtheorem{lemma}{Lemma}[section]

\numberwithin{equation}{section}

\pgfplotsset{compat=1.18}

\begin{document}
\begingroup
\def\uppercasenonmath#1{} 
\title{\large Inverse scattering for the multipoint potentials of Bethe$-$Peierls$-$Thomas$-$Fermi type}
\let\MakeUppercase\relax 

\author{\large \textit{P.C. Kuo},\, \textit{R.G. Novikov}}


\begin{abstract}
We consider the Schrödinger equation with a multipoint potential of the Bethe$-$Peierls$-$Thomas$-$Fermi type. We show that such a potential in dimension $d=2$ or $d=3$ is uniquely determined by its scattering amplitude at a fixed positive energy. Moreover, we show that there is no non-zero potential of this type with zero scattering amplitude at a fixed positive energy and a fixed incident direction. Nevertheless, we also show that a multipoint potential of this type is not uniquely determined by its scattering amplitude at a positive energy $E$ and a fixed incident direction. Our proofs also contribute to the theory of inverse source problem for the Helmholtz equation with multipoint source.
\end{abstract}

\maketitle


\section{Introduction}
\label{intro}
We consider the stationary Schrödinger equation in dimensions $d=1,2,3$,

\begin{equation}
    -\Delta \psi+\nu(x) \psi=E \psi, \quad x \in \mathbb{R}^{d}, \quad E>0,
    \label{Schrodinger}
\end{equation}
with a multipoint potential (scatterer) of Bethe$-$Peierls$-$Thomas$-$Fermi type:
 \begin{equation}
     \quad \nu(x)=\sum\limits_{j=1}^{n} \delta_{\alpha_{j}}\left(x-y_{j}\right), \quad \alpha_{j} \in \mathbb{C}, \quad y_{j} \in \mathbb{R}^{d}, \quad y_{j} \neq y_{j'} \quad \forall j \neq j'.
     \label{multipotential}
 \end{equation}
It is known that point scatterers $\delta_{\alpha}(x)$ of this type are only defined for $x \in \mathbb{R}^{d}$ with $d=1,2,3$. For $d=1$, $\delta_{\alpha}(x)=\epsilon \delta(x)$, where $\epsilon=-1/\alpha$ and $\delta$ denotes the standard Dirac delta function. For $d=2$ and $d=3$, the precise definition of "renormalized" delta function $\delta_{\alpha}$ is more subtle and is discussed, in particular, in \cite{albeverio2012solvable}, \cite{demkov2013zero}, \cite{grinevich2021transmission}, and in references therein. Historically, the aforementioned point scatterers in dimension $d=3$ were first introduced to describe the interaction between neutrons and protons by Bethe, Peierls \cite{bethe1935quantum}, Thomas \cite{thomas1935interaction} and Fermi \cite{fermi1936sul}. The mathematical theory of these point scatterers in dimensions $d=2$ and $d=3$ goes back to Zeldovich \cite{zel1960scattering}, Berezin and Faddeev\cite{Berezin:1960df}. Similar point scatterers also arise, in particular, in acoustics; see \cite{badalyan2009scattering}, \cite{dmitriev2021features} and references therein.

The simplest way to define a local solution $\psi$ of equation (\ref{Schrodinger}) with potential $\nu$ in (\ref{multipotential}) is as follows (see, for example, \cite{grinevich2021transmission}):

     \begin{equation}
         -\Delta \psi(x)=E \psi(x), \quad \forall x \in \mathbb{R}^{d} \backslash \{y_1,\cdots y_n\},
     \end{equation}
and, in addition:
\newline
if $d=1$, then $\psi(x)$ is continuous at $x=y_{j}$, but its first derivative has a jump
\begin{equation}
-\alpha_{j}\left[\psi^{\prime}\left(y_{j}+0\right)-\psi^{\prime}\left(y_{j}-0\right)\right]=\psi\left(y_{j}\right); 
\end{equation}
if $d=2$, then
\begin{equation}
\psi(x)=\psi_{j,-1} \ln \left|x-y_{j}\right|+\psi_{j, 0}+O\left(\left|x-y_{j}\right|\right) \quad \text {as} \quad x \rightarrow y_{j},
\end{equation}
with 
\begin{equation*}
 \psi_{j, 0}=(-2 \pi \alpha_{j}-\ln 2+\gamma);
\end{equation*}
if $d=3$, then
\begin{equation}
\psi(x)=\frac{\psi_{j,-1}}{\left|x-y_{j}\right|}+\psi_{j, 0}+O\left(\left|x-y_{j}\right|\right) \quad \text {as} \quad x \rightarrow y_{j},
\end{equation}
with 
\begin{equation*}
 4 \pi \alpha_{j}\psi_{j, -1}=\psi_{j, 0}.
\end{equation*}

We consider scattering and inverse scattering for model (\ref{Schrodinger}), (\ref{multipotential}). The direct scattering and spectral theory for this model is well-developed, at least, for real parameters $\alpha_{j}$'s in  (\ref{multipotential}); see, e.g., \cite{albeverio2012solvable}, \cite{Berezin:1960df}, \cite{demkov2013zero}, \cite{grinevich2024transparent}, \cite{grinevich2022spectral}, \cite{malamud2024kernels} and references therein. An important point is that the direct scattering problem for model (\ref{Schrodinger}), (\ref{multipotential}) is exactly solvable. Some formulas of this theory are recalled in Section \ref{Preliminaries}. These include formulas for the scattering amplitude $f^{+}$; see formulas (\ref{f+}), (\ref{scattering amplitude})-(\ref{coefficientA}).

For other exactly solvable models in quantum mechanics,
see. e.g., \cite{grinevich2000scattering},
\cite{grinevich2022spectral}, \cite{landau1958course}, \cite{novikov1984theory}, \cite{taimanov2010moutard} and references therein. 

In the present work we continue studies on inverse scattering for model (\ref{Schrodinger}), (\ref{multipotential}) in dimensions $d=2$ and $d=3$. 
In connection with results given in the literature in this direction, see \cite{agaltsov2019examples}, \cite{badalyan2009scattering}, \cite{dmitriev2021features}, \cite{gesztesy1999inverse}, \cite{grinevich2012faddeev}, \cite{grinevich2021transmission}, \cite{mantile2023inverse}, \cite{novikov2018inverse} and references therein.

We show that for model (\ref{Schrodinger}), (\ref{multipotential}), the scattering amplitude $f^{+}$ at a fixed positive energy $E$ uniquely determines the potential $\nu$ in (\ref{multipotential}); see Theorem \ref{inverse scattering} in Section \ref{Main}. Moreover, we show that there is no non-zero potential $\nu$ in (\ref{multipotential}) with zero scattering amplitude $f^{+}$ at a fixed positive energy $E$ and a fixed incident direction; see Theorem \ref{transparent} in Section \ref{Main}. Nevertheless, we also show that a multipoint potential $\nu$ in (\ref{multipotential}) is not uniquely determined by its scattering amplitude $f^{+}$ at a positive energy $E$ and a fixed incident direction; see Theorem \ref{counter} in Section \ref{Main}. In particular, we continue recent studies of \cite{grinevich2024transparent},
\cite{grinevich2021transmission} on possibilities for transparency or partial transparency for potential $\nu$ in (\ref{multipotential}) (that is vanishing or partial vanishing the scattering amplitude $f^{+}$). 

In addition, in connection with inverse scattering for the Schrödinger equation (\ref{Schrodinger}) with regular potential $\nu$, see  \cite{bukhgeim2008recovering}, \cite{chadan2012inverse}, \cite{grinevich2024transparent}, \cite{loran2024can}, \cite{novikov1988multidimensional}, \cite{novikov2019multidimensional}, and references therein. Note that studies on inverse scattering at fixed energy, including transparency, for the case of regular $\nu$ go back to Regge \cite{regge1959introduction}, R.Newton \cite{newton1962construction},  and Sabatier \cite{sabatier1966asymptotic}.

 Note that proofs of the present work also contribute to the theory of inverse source problem for the Helmholtz equation:

\begin{equation}
    \Delta \psi(x) +E \psi(x)=S(x), \quad x \in \mathbb{R}^{d}, \quad d \geq 2, \quad E>0,
    \label{Hemholtz}
\end{equation}
at fixed $E$, with $S$ of the form
\begin{equation}
    S(x)=\sum_{j=1}^{n} c_{j} \delta (x-y_{j}), \quad c_{j} \in \mathbb{C}, \quad y_{j} \in \mathbb{R}^{d},
    \label{Helmholtz source}
\end{equation}
where $\delta$ is the standard Dirac function. The related inverse source problem from far-field data consists in finding $S$ from the far-field amplitude $a^{+}$ for model (\ref{Hemholtz}), (\ref{Helmholtz source}) at fixed $E$.
Some formulas for $a$ are recalled in Section \ref{Preliminaries}; see formulas (\ref{eigenfunctions far field}), (\ref{scattering amplitude far field}). Our results on the aforementioned inverse source problem consist in Proposition \ref{key lemma complex}
in Section \ref{Main} and its proof in Section \ref{proof1} via some complex analysis.

In connection with results given in the literature on the inverse source problem from near-field data for model (\ref{Hemholtz}), (\ref{Helmholtz source}), see, for example, \cite{bao2021recovering},  \cite{el2012holder}, \cite{el2011inverse} and references therein. In addition, the inverse source problem for the simplest discrete Helmholtz equation is recently studied in \cite{novikov2024inverse}.

Further structure of the paper is as follows. Some preliminaries are recalled in Section \ref{Preliminaries}. The main results of the present work are stated in Section \ref{Main}. The proofs are given in Sections \ref{proof1} and \ref{proofs}.

\section{Preliminaries} \label{Preliminaries}

Let 
\begin{equation}
    \mathbb{S}_{r}^{d-1}=\{x \in \mathbb{R}^{d};|x|=r\},\quad r>0; \quad 
    \mathbb{S}^{d-1}=\mathbb{S}_{1}^{d-1}.
\end{equation}
For equation (\ref{Schrodinger}) with potential $\nu$ in (\ref{multipotential}), we consider the scattering eigenfunctions $\psi^{+}$ such that
\begin{equation} 
\label{wave}
\psi^{+}(x,k)=\psi_{0} +\psi^{sc}(x,k),\quad x \in \mathbb{R}^{d},\quad k \in \mathbb{S}_{\kappa}^{d-1}
\end{equation}
where $\kappa=\sqrt{E}$, $\psi_{0}(x,k)=e^{i k \cdot x}$ and $\psi^{sc}(x,k)$ satisfies the Sommerfeld radiation condition:
\begin{equation}
\label{Sommerfeld}
|x|^{\frac{d-1}{2}} \left( \frac{\partial}{\partial |x|} - i\kappa\right) \psi^{sc}(x,k) \xrightarrow[|x| \to \infty]{} 0 \quad \text{uniformly in all directions } \frac{x}{|x|}.
\end{equation}
Then,
\begin{equation} 
\label{f+}
\psi^{sc}(x,k)= \frac{e^{i\kappa|x|}}{|x|^{(d-1) / 2}} f^{+}(k, \frac{\kappa}{|x|}x)+O\left(\frac{1}{|x|^{(d+1) / 2}}\right) \quad \text{as } \quad |x| \rightarrow \infty.
\end{equation}
The function $f^{+}$ arising in (\ref{f+}) is defined on $\mathbb{S}_{\kappa}^{d-1}\times\mathbb{S}_{\kappa}^{d-1}$, and is the scattering amplitude for model (\ref{Schrodinger}),(\ref{multipotential}). It is convenient to present $f^{+}$ as follows:

\begin{equation}
f^{+}(k,l)=c(d,\kappa)f(k,l),\quad k,l \in \mathbb{S}^{d-1}_{\kappa}; \quad c(d,\kappa) =-\pi i(\sqrt{2 \pi} e^{-i \pi / 4} )^{(d-1)}\kappa^{(d-3) / 2}.  
\end{equation}

Crucially, when $\nu$ is the multipoint potential of Bethe–Peierls–Thomas$-$Fermi type (\ref{multipotential}), we have the following explicit formulas for the scattering eigenfunction $\psi^{+}$ and scattering amplitude $f^{+}$ (see \cite{grinevich2021transmission} and references therein):

\begin{align}
\psi^{+}(x, k) &=e^{i k \cdot x}+\sum_{j=1}^{n} q_{j}(k) G^{+}\left(x-y_{j},\kappa\right), \quad x \in \mathbb{R}^{d},\quad k \in \mathbb{S}^{d-1}_{\kappa}, \label{eigenfunctions} \\
f(k, l) &=\frac{1}{(2 \pi)^{d}} \sum_{j=1}^{n} q_{j}(k) e^{-i l \cdot y_{j}},\quad k,l \in \mathbb{S}^{d-1}_{\kappa},  \label{scattering amplitude}
\end{align}
where  $q(k)=\left(q_{1}(k), \cdots, q_{n}(k)\right)^{t}$ satisfies
\begin{align}
    A(\kappa) q(k) &=b(k), \label{A=qb} \\
    A_{j, j}  (\kappa) &= 
    \begin{cases}\alpha_{j}+(2 i \kappa)^{-1} &, \text {if } d=1 \\ \alpha_{j}-(4 \pi)^{-1}(\pi i-2 \ln (\kappa)) &, \text {if } d=2 \\
    \alpha_{j}-i(4 \pi)^{-1}\kappa &, \text {if } d=3 \\
    \end{cases} \label{coefficientA} \\
    A_{j, j^{\prime}}(\kappa) &= G^{+}\left(y_{j}-y_{j^{\prime}}, \kappa \right), j \neq j^{\prime}, \nonumber \\
    b(k) &=-\left(e^{i k \cdot y_{1}}, e^{i k \cdot y_{2}}, \cdots e^{i k \cdot y_{n}}\right)^{t}, \nonumber
\end{align}
and $G^{+}$ is the  Green function with the Sommerfeld radiation condition for the operator $\Delta+E$, given by

\begin{align}
G^{+}(x, \kappa)&=-(2 \pi)^{-d} \int_{\mathbb{R}^{d}} \frac{e^{i \xi \cdot x} d \xi}{|\xi|^{2}-\kappa^2-i \cdot 0}    \label{Green} \\
&= 
\begin{cases}
  \frac{e^{i\kappa|x|}}{2 i\kappa}  &  ,\text{if } d=1  \\
  -\frac{i}{4} H_{0}^{1}(|x|\kappa) & ,\text{if }d=2 \\
  -\frac{e^{i\kappa|x|}}{4 \pi|x|} & ,\text{if }d=3 \\ 
\end{cases} \nonumber
\end{align}
where $H_{0}^{1}$ denotes the zeroth-order Hankel function of the first kind.
Note that $G^{+}(x,\kappa)$ only depends on $|x|$ when $E=\kappa^2$ is fixed.

In connection with singularity, the following formulas hold:

\begin{align*}
 -\frac{1}{4 \pi} q_{j}&=\quad \psi_{j,-1}, \quad d=3 \\
 \frac{1}{2 \pi} q_{j}&=\quad \psi_{j,-1}, \quad d=2 \\
 q_{j}&=\quad \psi^{\prime}\left(y_{j}+0\right)-\psi^{\prime}\left(y_{j}-0\right), \quad d=1
\end{align*}

Note that we assume the condition
\begin{equation}
\label{det}
    \operatorname{det}(A(\kappa)) \neq 0,
\end{equation}
which is automatically true when $\alpha_{j}$'s are all real.

For equation (\ref{Hemholtz}) with the source function $S$ in (\ref{Helmholtz source}), we consider the solution $\psi$ satisfying the Sommerfeld radiation condition (\ref{Sommerfeld}), where $\kappa=\sqrt{E}>0$. Then

\begin{equation}
\label{helmholtz scattering amplitude}
    \psi(x,\kappa)= \frac{e^{i\kappa|x|}}{|x|^{(d-1) / 2}} a^{+}(\frac{\kappa}{|x|}x)+O\left(\frac{1}{|x|^{(d+1) / 2}}\right) \quad \text{as } \quad |x| \rightarrow \infty.
\end{equation}

The function $a^{+}$ arising in ($\ref{helmholtz scattering amplitude}$) is defined on $\mathbb{S}_{\kappa}^{d-1}$ and is
the far-field amplitude for model (\ref{Hemholtz}), (\ref{Helmholtz source}).
In a similar way with $f^{+}$, it is convenient to present $a^{+}$ as follows:

\begin{equation}
\label{a+}
a^{+}(l)=c(d,\kappa)a(l),\quad l \in \mathbb{S}^{d-1}_{\kappa}; \quad c(d,\kappa) =-\pi i(\sqrt{2 \pi} e^{-i \pi / 4} )^{(d-1)}\kappa^{(d-3)/2}.   
\end{equation}

We have the following formulas for
$\psi$ and $a$:

\begin{align}
\psi(x, \kappa) &=\sum_{j=1}^{n} c_{j} G^{+}\left(x-y_{j},\kappa\right), \quad x \in \mathbb{R}^{d} \label{eigenfunctions far field}, \\
a(l) &=\frac{1}{(2 \pi)^{d}} \sum_{j=1}^{n} c_{j} e^{-i l \cdot y_{j}},\quad l \in \mathbb{S}^{d-1}_{\kappa}. \label{scattering amplitude far field}
\end{align}

One can see a similarity between formulas (\ref{eigenfunctions}), (\ref{scattering amplitude}) for model (\ref{Schrodinger}), (\ref{multipotential}) and formulas (\ref{eigenfunctions far field}), (\ref{scattering amplitude far field}) for model (\ref{Hemholtz}), (\ref{Helmholtz source}).

\section{Main results}\label{Main}

We start with the following global uniqueness result on inverse scattering 
for model (\ref{Schrodinger}), (\ref{multipotential}) at a fixed energy $E>0$ in dimension $d=2$ or $d=3$.
\begin{theorem}
\label{inverse scattering}
    A multipoint potential $\nu$ of the form (\ref{multipotential}), for $d=2$ or $d=3$ under condition (\ref{det}), is uniquely determined by its scattering amplitude $f^{+}$ at a fixed energy $E>0$.
\end{theorem}

Prototypes of Theorem \ref{inverse scattering} for regular potentials $\nu$ are given in \cite{bukhgeim2008recovering} and \cite{novikov2019multidimensional}. This theorem is proved in Section \ref{proofs}. In turn, this proof uses the following result of independent interest.

\begin{proposition}
\label{key lemma complex}
Let $y_{1}, y_{2}, \cdots, y_{n} \in$ $\mathbb{R}^{d}$, with $d \geqslant 2$, be mutually distinct and $c_{1}, c_{2}, \cdots, c_{n} \in \mathbb{C}\backslash \left\{0\right\} $, where $n \in \mathbb{N}$. Let
$u(\theta)=\sum\limits_{j=1}^{n} c_{j} e^{i y_{j} \cdot \theta}, \theta \in \mathbb{S}^{d-1}$. Then
$y_j$'s and $c_j$'s are uniquely determined by $u$ on $\Theta$, where $\Theta$ is a non-empty open subset of $ \mathbb{S}^{d-1}$.
\end{proposition}

Proposition \ref{key lemma complex} is proved in Section \ref{proof1} using some complex analysis. Of course, this proposition for real $y_{j}$'s can also be proved by using relations between far-field and near-field data and results of \cite{bao2021recovering}, \cite{fermi1936sul}, \cite{thomas1935interaction}.
However, the later proof looks much more complicated. 

Note that Proposition \ref{key lemma complex} also holds for $y_{j} \in \mathbb{C}^{d}$. 

In virtue of (\ref{helmholtz scattering amplitude}), (\ref{a+}), (\ref{scattering amplitude far field}), Proposition \ref{key lemma complex} implies the global uniqueness for the inverse source problem from far-field data for model (\ref{Hemholtz}), (\ref{Helmholtz source}) at fixed $E$. Note that Proposition \ref{key lemma complex} is not valid even for real $y_{j}$'s for some other convex and real-analytic surfaces in place of $\mathbb{S}^{d-1}$. In particular, article \cite{novikov2024inverse} gives related examples of non-uniqueness in the case of the surfaces
$$\Gamma(\lambda)=\{k \in [-\pi,\pi]^{d} ;\sum\limits_{i=1}^{n}2\cos(k_{i})=\lambda  \}, \text{ where } |\lambda| \in (2d-4,2d).$$

Due to Theorem \ref{inverse scattering}, there are no non-zero potentials $\nu$ in (\ref{multipotential}) which are 
transparent at a fixed positive energy $E$ in dimensions $d = 2, 3$ under assumption (\ref{det}), that is with $f^+ \equiv 0$ at fixed $E$. Moreover, the following much stronger result in this direction also holds.

\begin{theorem}
 Let $\nu$ be a multipoint potential of the form (\ref{multipotential}), with $d=2$ or $d=3$  under condition (\ref{det}), and $f^{+}$ be its scattering amplitude. Suppose that $f^{+}(k,l) \equiv 0$ at a fixed energy $E>0$ and a fixed incident vector $k$, then $\nu \equiv 0$.
 \label{transparent}
\end{theorem}

In the framework of formula (\ref{multipotential}) and
related formulas in Sections \ref{intro} and \ref{Preliminaries},
the case of $\nu \equiv 0$ can be considered as 
the case of $n=0$ or, in other words,
as the case of no scatterers.
The case of $\nu \equiv 0$ can also be considered as the case when $\alpha_j = \infty$ for all $j$.

In spite of Theorem \ref{transparent}, the following non-uniqueness result also holds.

\begin{theorem}

\label{counter}
    A multipoint potential $\nu$ of the form (\ref{multipotential}), for $d=2$ or $d=3$ under condition (\ref{det}), is not uniquely determined, in general, by its scattering amplitude $f^{+}(k,l)$ at a positive energy $E$ and a fixed incident vector $k$.
\end{theorem}
Theorem \ref{transparent} and Theorem \ref{counter} are proved in Section \ref{proofs}.

In particular,  by the results of Theorem \ref{inverse scattering}-\ref{counter} we continue recent studies of \cite{grinevich2024transparent}, \cite{grinevich2021transmission}, 
on possibilities of transparency or partial transparency for potential $\nu$ in (\ref{multipotential}).

\section{Proof of Proposition \ref{key lemma complex}} \label{proof1}

In addition to $u$ on $\mathbb{S}^{d-1}$, we also consider
\begin{equation}    u(\theta)=\sum\limits_{j=1}^{n}c_{j}e^{iy_{j}\cdot \theta},\quad \theta \in \Sigma^{d-1},
\end{equation}
where
\begin{equation}
    \Sigma^{d-1}=\left\{\theta \in \mathbb{C}^{d} ; \theta^{2}=\theta_{1}^{2}+\cdots+\theta_{d}^{2}=1\right\}.
\end{equation}

Note that $\Sigma^{d-1} \cap \mathbb{R}^{d}=\mathbb{S}^{d-1}$. 

As $u$ is holomorphic on $\Sigma^{d-1}$, the restriction $u|_{\Theta}$ uniquely determines $u$ on ${\mathbb{S}^{d-1}}$ and on ${\Sigma^{d-1}}$ via analytic continuation. 

Consider the complex vector-valued function
\begin{equation}
\label{theta}
    \theta(\tau,e_{1},e_{2})= \sqrt{1+\tau^{2}}e_{1}-i\tau e_{2}, \quad \tau \in \mathbb{R}, \quad (e_{1},e_{2}) \in U\mathbb{S}^{d-1},
\end{equation}

where $U\mathbb{S}^{d-1}=\{(e_{1},e_{2})\in \mathbb{S}^{d-1}\times\mathbb{S}^{d-1}; e_{1}\cdot e_{2}=0\}$ denotes the unit tangent bundle on $\mathbb{S}^{d-1}$. 

One can see that $\theta$ defined in (\ref{theta}) takes values in $\Sigma^{d-1}$.

Then
\begin{equation}
\label{utheta}
u(\theta(\tau,e_{1},e_{2}))=\sum_{i=1}^{n}c_{j}e^{iy_{j}\cdot \theta(\tau,e_{1},e_{2})}=\sum_{i=1}^{n}c_{j}e^{\tau e_{2}\cdot y_{j}} e^{-i \sqrt{1+\tau^{2}}e_{1}\cdot y_{j}}, \quad \tau \in \mathbb{R}, \quad (e_{1},e_{2}) \in U\mathbb{S}^{d-1},
\end{equation}
where the modulus of each term is
\begin{equation}
    |c_{j}e^{iy_{j}\cdot \theta(\tau,e_{1},e_{2})}|= |c_{j}| e^{\tau e_{2}\cdot y_{j}}.
\end{equation}

Now we propose the following recursive determination of all $y_{j}$'s and $c_{j}$'s. The recursive determination is based on the behavior of $\tilde{u}:=u \circ \theta$ on $\mathbb{R}\times U\mathbb{S}^{d-1}$, which splits into the following three cases:

\begin{enumerate}[(i)]
    \item $\tilde{u} \equiv 0$. In view of (\ref{utheta}), this implies $n=0$ and $u \equiv 0$. The determination is completed.
    \item $\tilde{u} \not\equiv 0$ but $\tilde{u}$ is bounded on $\mathbb{R}\times U\mathbb{S}^{d+1}$. In view of (\ref{utheta}), this implies $n=1$, $y_{1}=0$, and $\tilde{u} \equiv c_{1}$. The determination is completed. 
    \item $\tilde{u}$ is unbounded. In view of (\ref{utheta}), we can
    take $e_{2} \in \mathbb{S}^{d-1}$ such that $\tilde{u}$ has maximal exponential growth as $\tau$ tends to $+\infty$. This implies that $e_{2}=\frac{y_{l}}{|y_{l}|}$ for some $y_{l}$ such that $|y_{l}|=\max\limits_{1\leqslant k \leqslant n} |y_{k}|\neq 0$. In this way, for any fixed $e_{1}$ such that $(e_{1},e_{2})\in U\mathbb{S}^{d-1}$, we have $e_{1}\cdot y_{l}=0$ and, in addition, $e_{2}\cdot y_{l}=|y_{l}| > |e_{2}\cdot y_{k}|$ for all $k  \neq l$. Therefore,
    \begin{equation}
    \label{asymu}
    \tilde{u}(\tau,e_{1},e_{2}) \sim c_{l}e^{\tau |y_{l}|} \quad \text{as } \tau \to \infty.  
    \end{equation}
Now we reconstruct $y_{l}$ and $c_{l}$ using (\ref{asymu}) and using that $e_{2}=\frac{y_{l}}{|y_{l}|}$. We complete case (iii)  by setting
$$
u_{new}(\theta):=u(\theta)-c_{l}e^{iy_{l}\cdot\theta}=\sum\limits_{1\leqslant j \leqslant n, j \neq l}c_{j}e^{iy_{j}\cdot\theta},\quad \theta \in \Sigma^{d-1}.
$$ 
Then we proceed in the same way based on the behavior of $\tilde{u}_{new}=u_{new}\circ \theta$.
\end{enumerate}
This recursive determination will terminate within $n$ or $n+1$ iterations and gives all $y_{j}$'s and $c_{j}$'s.

\section{Proofs of Theorem \ref{inverse scattering}, \ref{transparent}}
\label{proofs}

\subsection{Nontriviality of  $q_{j}$}

We start with the following lemma.
\begin{lemma}
\label{q_j(k)}
Under the assumption of Theorem \ref{inverse scattering}, for each $j=1,\cdots,n$, the coefficient $q_{j}$ in formula (\ref{scattering amplitude}) is an analytic function on $\mathbb{S}_{\kappa}^{d-1}$ and $q_{j}$ is not identically zero.
\end{lemma}

By relation (\ref{A=qb}) and condition (\ref{det}), for all $1\leqslant j \leqslant n$, we have
\begin{equation}
\label{qjsum}
q_{j}(k)=\sum\limits_{l=1}^{n} (A(\kappa))^{-1}_{j,l} e^{ik\cdot y_{l}}, \quad k \in \mathbb{S}_{\kappa}^{d-1}.
\end{equation}

For a fixed $j$, not all $(A(\kappa))^{-1}_{j,1},\cdots,(A(\kappa))^{-1}_{j,n}$ are zero since $\det A^{-1}(\kappa) = 
(\det A(\kappa))^{-1} \neq 0$. Therefore, Lemma \ref{q_j(k)} follows from formula (\ref{qjsum}) and from Proposition \ref{key lemma complex} for $u(\theta)=q_{j}(\kappa \theta)$.

Note that in view of analycity of $q_j$, we also have that $q_{j}(k) \neq 0$ almost everywhere on $S_{\kappa}^{d-1}$.

\subsection{Proof of Theorem \ref{inverse scattering}}

Applying Proposition \ref{key lemma complex} and Lemma \ref{q_j(k)} to $u(\theta)=f(k,\kappa \theta)$ given by formula (\ref{scattering amplitude}), one can see that $f$ at fixed $E>0$ uniquely determines $q_{j}$ on $\mathbb{S}^{d-1}_{\kappa}$ and $y_{j}$ for all $1\leqslant j \leqslant n$.

Next, when $q_{j}(k) \neq 0$, the diagonal entry $A_{j,j}(\kappa)$ is determined by the 
\begin{equation}
\label{Ajj}
    A_{j,j}(\kappa)q_{j}(k)=-e^{ik\cdot y_{j}}-\sum\limits_{j' \neq j} A_{j,j'}(\kappa) q_{j'}(k), \quad 1 \leqslant j \leqslant n.
\end{equation}
Here, equation (\ref{Ajj}) follows from formula (\ref{A=qb}).
Finally, $A_{j,j}(\kappa)$ determines  $\alpha_{j}$ via formula (\ref{coefficientA}) for all $1 \leqslant j \leqslant n$.
This completes the proof of Theorem \ref{inverse scattering}.

\subsection{Proof of Theorem \ref{transparent}}

 Let $E>0$ and $k \in \mathbb{S}^{d-1}_{\kappa}$. Suppose that $f^{+}(k,l) = 0$ for all $l \in \mathbb{S}^{d-1}_{\kappa}$. By formula (\ref{scattering amplitude}) and Proposition \ref{key lemma complex}, $q_{j}(k)=0$ for all $1\leqslant j \leqslant n$. It follows from relation (\ref{A=qb}) that $b(k)=0$.
 However, $b_{j}(k)=-e^{ik\cdot y_{j}} \neq 0$ for all $1\leqslant j \leqslant n$. This implies that $v \equiv 0$.

\section{Proofs of Theorem \ref{counter}}

In this section, we provide two approaches to constructing counterexamples. The first one is based on fitting parameters of two-point scatterers. The second one consists of adding a point scatterer at a special position which corresponds to a zero of the wave function $\psi^{+}$ in (\ref{eigenfunctions}).
\subsection{Proof of Theorem \ref{counter} by fitting parameters.}
\label{fitting}
Consider $\nu$ and $\tilde{\nu}$ with parameters 
$n$, $\alpha_j$, $y_j$ and $\tilde{n}$, $\tilde{\alpha}_j$, $\tilde{y}_j$
respectively, where $n=\tilde{n}= 2$, $y_{1}=\tilde{y}_{1}=0$, $y_{2}=-\tilde{y}_{2}\neq 0$, $ky_{2}=0$.

\begin{equation*}
\alpha_{1}=\tilde{\alpha}_{1}=
\begin{cases}
  \frac{1}{4\pi}(\pi i -2 \ln(\kappa))-\frac{i}{4} H_{0}^{1}(|y_{2}|\kappa) & ,d=2, \\
  \frac{i\kappa}{4 \pi}-\frac{e^{i\kappa|y_{2}|}}{4\pi |y_{2}|} & ,d=3, \\
\end{cases}
\end{equation*}
and $\alpha_2$ and $\tilde{\alpha}_2$ are arbitrarily complex numbers other than $\alpha_1$. By (\ref{coefficientA}), (\ref{Green}) and our parameters chosen, we have 
$$
A=
\begin{pmatrix}
g & g \\
g & h 
\end{pmatrix},\quad 
\tilde{A}=
\begin{pmatrix}
g & g \\
g & \tilde{h}
\end{pmatrix}, \text{ and }
b=\tilde{b}=
\begin{pmatrix}
-1 \\
-1 
\end{pmatrix},
$$
where $g=G^{+}(y_{2},\kappa) \neq 0$, $h \neq g$, and $\tilde{h} \neq g$, ensuring that  $\det A \neq 0$ and $\det \tilde{A} \neq 0$. 

Then equation (\ref{A=qb}) implies
$$
q(k)=\tilde{q}(k)=\frac{-1}{g}
\begin{pmatrix}
1 \\
0 
\end{pmatrix}.
$$
Therefore, by formula (\ref{scattering amplitude}),
$$
f^{+}(k,l)=\frac{-1}{g (2 \pi)^{d}}=\tilde{f}^{+}(k,l), \quad \forall l \in \mathbb{S}^{d-1}_{\kappa}.
$$
Nevertheless, $v \not \equiv \tilde{v}$ as $y_{2} \neq\tilde{y}_{2}$.

Thus, Theorem \ref{counter} is proved.
However, by this fitting approach, we did not find counterexamples when all $\alpha_{j}$'s and $\tilde{\alpha}_{j}$'s are real.

\subsection{Proof of Theorem \ref{counter} by adding invisible point scatterers} \label{adding invisible}
We use the notations $\nu$, $n$, $y_{j}$, $\alpha_{j}$ and $\tilde{\nu}$, $\tilde{n}$, $\tilde{y}_{j}$, $\tilde{\alpha}_{j}$ in a similar way with Subsection \ref{fitting}.
Let $\nu$ be a multipoint potential of (\ref{multipotential}). Suppose there exists $\tilde{y}_{n+1}\in \mathbb{R}^d \backslash \{ y_{1},...,y_{n} \}$ such that
$\psi^{+}(y_{n+1}, k)=0$, where $\psi^{+}$ is given by formula (\ref{eigenfunctions}); that is,

\begin{equation}
    e^{i k \tilde{y}_{n+1}}+\sum_{j=1}^{n} q_{j}(k) G^{+}\left(\tilde{y}_{n+1}-y_{j},\kappa\right)=0,\quad\tilde{y}_{n+1}\in \mathbb{R}^d \backslash \{ y_{1},...,y_{n} \}.
    \label{invisible}
\end{equation}
Let $y_{j}=\tilde{y}_{j},\,1 \leqslant j \leqslant n$. Next, choose $\tilde{\alpha}_{n+1} \in \mathbb{C}$ so that $\det{\tilde{A}} \neq 0$. This is possible since $\det{A} \neq 0$ and by formula (\ref{coefficientA}), we have 

\begin{equation}
\label{Atilde}
\tilde{A}=
\begin{pmatrix} 
  A & g^{t}  \\ 
  g  & \tilde{A}_{n+1,n+1}
\end{pmatrix},
\end{equation}
where 
$$
\tilde{A}_{n+1,n+1}=
\begin{cases} \tilde{\alpha}_{n+1}-(4 \pi)^{-1}(\pi i-2 \ln (\kappa)) &, \text {if } d=2, \\
\tilde{\alpha}_{n+1}-i(4 \pi)^{-1}\kappa &, \text {if } d=3,
\end{cases}
$$
and
$$
g=\left( G^{+}\left(\tilde{y}_{1}-y_{n+1},\kappa\right),\cdots,G^{+}\left(\tilde{y}_{n}-y_{n+1},\kappa\right) \right).
$$
It follows from (\ref{A=qb}), (\ref{invisible}), and (\ref{Atilde}) that
$$
\tilde{q}=
\begin{pmatrix} 
  q  \\ 
  0 
\end{pmatrix}.
$$

Therefore, by (\ref{scattering amplitude}), we have $\tilde{f^{+}}(k,l)=f^{+}(k,l)$ for all $l \in \mathbb{S}^{d-1}_{\kappa}$. This completes the schema of proof of
Theorem \ref{counter} by adding an invisible point scatterer at a zero of $\psi^{+}$.

In addition, by (\ref{eigenfunctions}), we also have that $\tilde{\psi}^{+}(x,k)=\psi^{+}(x,k)$ for all $x \in \mathbb{R}^{d}$. Therefore, we can continue adding invisible point scatterers if $\psi^{+}$
has other zeros; i.e., equation (\ref{invisible}) has other solutions.

In general, $\tilde{y}_{n+1}$ satisfying (\ref{invisible}) may not exist. However, there are already the following examples with $n=1$ when such $\tilde{y}_{n+1}$ exists. 

Let $d=2$, $n=1$, $y_1=0$, and $\alpha_1=\frac{\pi i-2 \ln (\kappa)}{4 \pi}-\frac{iH_{0}^{1}(1)}{4e^{i}}$. We have $A=\alpha_1-\frac{\pi i-2 \ln (\kappa)}{4 \pi}=-\frac{iH_{0}^{1}(1)}{4e^{i}}$ by formula (\ref{coefficientA}) and $q(k)=-A^{-1}e^{ik\dot y_1}=\frac{4e^{i}}{iH_{0}^{1}(1)}$ by formula (\ref{A=qb}). Now, for $\tilde{y}_{2}=\frac{k}{\kappa^2}$, we have that
$$
e^{ik\tilde{y}_{2}}+q(k)G^{+}(\tilde{y}_{2},\kappa)=e^{i}+ \frac{4e^{i}}{iH_{0}^{1}(1)} G^{+}(\frac{k}{\kappa^{2}},\kappa)=e^{i}+\frac{4e^{i}}{iH_{0}^{1}(1)} \frac{-i}{4}H_{0}^{1}(1)=0.
$$
Thus, this $\tilde{y}_{2}$ is a solution of (\ref{invisible}).

Let $d=3$, $n=1$, $y_1=0$, and $\alpha_1=\frac{i\kappa}{4\pi}-\frac{1}{4 \pi \kappa}$. We have $A=\alpha_1-\frac{i\kappa}{4 \pi}=-\frac{1}{4\pi\kappa}$ by formula (\ref{coefficientA}) and $q(k)=-A^{-1}e^{ik\dot y_1}=4\pi\kappa$ by formula (\ref{A=qb}). Now, for $\tilde{y}_{2}=k$, we have that
$$
e^{ik\tilde{y}_{2}}+q(k)G^{+}(\tilde{y}_{2},\kappa)=e^{i\kappa^{2}}+ \frac{1}{4\pi\kappa} G^{+}(k,\kappa)=e^{i\kappa^{2}}+4\pi\kappa\frac{-e^{i\kappa^{2}}}{4\pi\kappa}=0.
$$
Therefore, $\tilde{y}_{2}=k$ is a solution of (\ref{invisible}). 

However, in these examples, $\alpha_{1}$ is not real. 

In this connection, in Section \ref{numerical}, we give an example with real $\alpha_{j}$'s through numerical implementation.

\section{Numerical example}
\label{numerical}
This section aims to show solutions of equation (\ref{invisible}) by numerical simulation for the case when $\alpha_{j}$'s are all real. The purpose is to find a zero of the following function: 
\begin{equation}
   \psi^{+}(x):= e^{i k \cdot x}+\sum_{j=1}^{n} q_{j}(k) G^{+}\left(x-y_{j},\kappa\right), \quad x \in \mathbb{R}^d.  
\end{equation}
Here, we are interested in the case when $\alpha$'s are real since we have seen examples with complex $\alpha$'s.

Let $d=2$, $n=2$, $\alpha_{1}=1$, $\alpha_{2}=1$, $ k=\begin{pmatrix}
    0.05 \\ 
    0.05
\end{pmatrix}$,  $y_{1}=\begin{pmatrix} 
  0  \\ 
  0 
\end{pmatrix}$, and $y_{2}=\begin{pmatrix} 
  0.1  \\ 
  0.15 
\end{pmatrix}$. In this case, $q = \begin{pmatrix}
    0.51 - 1.86i
    \\
    0.51 - 1.86i
\end{pmatrix}$. In addition, we find a solution of (\ref{invisible}) around the point $(0.994,-4.398)$ by tracing the level lines of the real part and imaginary part of $\psi^{+}$. This is illustrated in the figure  \ref{reim}. This example show that the approach of adding invisible point scatterers in Subsection \ref{adding invisible} also works for the case when all $\alpha_j$'s are real.

\begin{figure}[hbt]
\caption{Level sets of $\Re (\psi^{+})=0$ and $\Im (\psi^{+})=0$ intersect at a point around $(0.994,-4.398)$}
\label{reim}
\includegraphics[width=8cm]{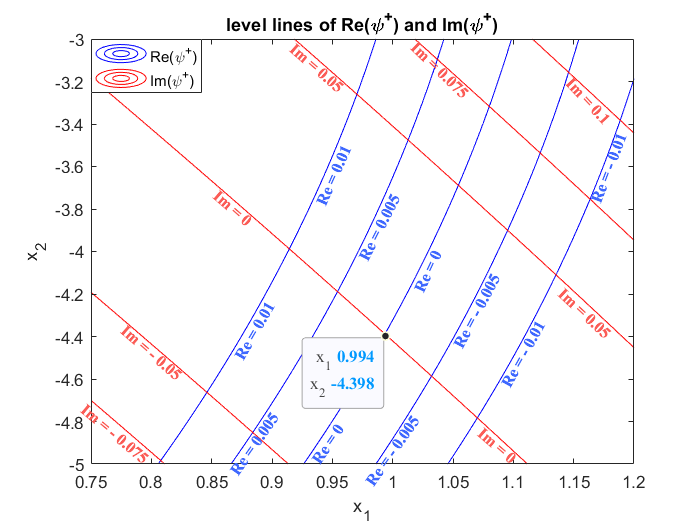}
\end{figure}
\newpage
\section*{Acknowledgements}

The main part of this work was fulfilled during the internship of the first author in the Centre de Mathématiques Appliquées of École polytechnique in June-December 2024 in the framework of research program for international talents.

\bibliographystyle{plain}     
\bibliography{scatter}   

\vspace{2cm}
\noindent
Pei-Cheng Kuo\\
Department of Mathematics, National Taiwan University,\\
4, Roosevelt Rd., Taipei 106319, Taiwan;\\
Université de Versailles Saint-Quentin-en-Yvelines, Université Paris-Saclay,\\
 45 Av. des États Unis, 78000 Versailles, France;\\
E-mail: r12221017@ntu.edu.tw
\\
\\

\noindent
Roman G. Novikov\\
CMAP, CNRS, École polytechnique,\\
Institut polytechnique de Paris, 91120 Palaiseau, France;\\
IEPT RAS, 117997 Moscow, Russia;\\
E-mail: novikov@cmap.polytechnique.fr
\end{document}